\documentclass[onefignum,onetabnum]{siamonline250211}

\usepackage[export]{adjustbox}

\usepackage{lipsum}
\usepackage{amsfonts}
\usepackage{graphicx}
\usepackage{epstopdf}
\ifpdf
  \DeclareGraphicsExtensions{.eps,.pdf,.png,.jpg}
\else
  \DeclareGraphicsExtensions{.eps}
\fi

\usepackage{enumitem}
\setlist[enumerate]{leftmargin=.5in}
\setlist[itemize]{leftmargin=.5in}

\newsiamremark{remark}{Remark}
\newsiamremark{hypothesis}{Hypothesis}
\crefname{hypothesis}{Hypothesis}{Hypotheses}
\newsiamthm{claim}{Claim}
\newsiamremark{fact}{Fact}
\crefname{fact}{Fact}{Facts}

\headers{Multiway Synergy Analysis}{E. Newman, A. Arnold, D. Calvetti, and E. Somersalo}

\title{Multiway Analysis and Uncertainty Quantification of Muscle Activation Patterns 
\thanks{
\funding{This material is based upon work supported by the National Science Foundation under Grant No. DMS-1929284 while the authors were in residence at the Institute for Computational and Experimental Research in Mathematics (ICERM) in Providence, RI, during the ``Stochastic and Randomized Algorithms in Scientific Computing: Foundations and Applications'' semester program. 
D. Calvetti was partly supported by the National Science Foundation (NSF) grant DMS-2513481, as well as by the Simons Foundation Fellowship. 
E. Somersalo was partly supported by the NSF grants DMS-2204618 and DMS-2513481.
E. Newman was partly support from the United States Air Force Office of Scientific Research under FA9550-26-1-B049 (program manager Dr. Fariba Fahroo) and by the NSF under grants DMS-2309751 and DMS-2541280. 
The support is acknowledged with gratitude.}}}

\author{Elizabeth Newman\thanks{Tufts University, Medford, MA 
  (\email{e.newman@tufts.edu}, \url{elizabethnewman.owlstown.com}).}
\and Andrea Arnold\thanks{Worcester Polytechnic Institute, Worcester, MA
  (\email{anarnold@wpi.edu}).}
\and Daniela Calvetti\thanks{Case Western Reserve University, Cleveland, OH 
  (\email{dxc57@case.edu}).}
\and Erkki Somersalo\footnotemark\thanks{Case Western Reserve University, Cleveland, OH 
  (\email{ejs49@case.edu}).}
  }

\usepackage{amsopn}

\DeclareMathOperator*{\fitmetric}{fit}

\usepackage{subcaption}
\usepackage{wrapfig}
\usepackage{stmaryrd} 
\usepackage{mathtools} 
\usepackage{bbm}
\usepackage{algorithm}
\usepackage{algpseudocode}
\usepackage{mathrsfs}
\usepackage{booktabs}
\usepackage[table]{xcolor}
\usepackage{xspace}

\usepackage{biblatex}
\def\mydefb#1{\expandafter\def\csname bf#1\endcsname{\mathbf{#1}}}
\def\mydefallb#1{\ifx#1\mydefallb\else\mydefb#1\expandafter\mydefallb\fi}
\mydefallb aAbBcCdDeEfFgGhHiIjJkKlLmMnNoOpPqQrRsStTuUvVwWxXyYzZ\mydefallb

\def\mydefb#1{\expandafter\def\csname #1bb\endcsname{\mathbb{#1}}}
\def\mydefallb#1{\ifx#1\mydefallb\else\mydefb#1\expandafter\mydefallb\fi}
\mydefallb aAbBcCdDeEfFgGhHiIjJkKlLmMnNoOpPqQrRsStTuUvVwWxXyYzZ\mydefallb

\def\mydefb#1{\expandafter\def\csname #1cal\endcsname{\mathcal{#1}}}
\def\mydefallb#1{\ifx#1\mydefallb\else\mydefb#1\expandafter\mydefallb\fi}
\mydefallb aAbBcCdDeEfFgGhHiIjJkKlLmMnNoOpPqQrRsStTuUvVwWxXyYzZ\mydefallb

\def\mydefgreek#1{\expandafter\def\csname bf#1\endcsname{\text{$\boldsymbol{\csname #1\endcsname}$}}}
\def\mydefallgreek#1{\ifx\mydefallgreek#1\else\mydefgreek{#1}%
   \lowercase{\mydefgreek{#1}}\expandafter\mydefallgreek\fi}
\mydefallgreek {alpha}{Alpha}{beta}{Beta}{gamma}{Gamma}{delta}{Delta}{epsilon}{Epsilon}{zeta}{Zeta}{eta}{Eta}{theta}{Theta}{iota}{Iota}{kappa}{Kappa}{lambda}{Lambda}{mu}{Mu}{nu}{Nu}{omicron}{Omicron}{pi}{Pi}{rho}{Rho}{sigma}{Sigma}{tau}{Tau}{upsilon}{Upsilon}{phi}{Phi}{xi}{Xi}{chi}{Chi}{psi}{Psi}{omega}{Omega}\mydefallgreek

\def\mydefb#1{\expandafter\def\csname T#1\endcsname{\boldsymbol{\mathcal{\MakeUppercase{#1}}}}}
\def\mydefallb#1{\ifx#1\mydefallb\else\mydefb#1\expandafter\mydefallb\fi}
\mydefallb aAbBcCdDeEfFgGhHiIjJkKlLmMnNoOpPqQrRsStTuUvVwWxXyYzZ\mydefallb

\def\ntime{{n_t}}
\def\nmuscle{{n_m}}
\def\nsample{{n_s}}

\def\Utime{\bfU_t}
\def\Umuscle{\bfU_m}
\def\Usample{\bfU_s}

\def\Wtime{\bfW_t}
\def\Wmuscle{\bfW_m}
\def\Wsample{\bfW_s}

\def\Xtime{\bfX_t}
\def\Xmuscle{\bfX_m}
\def\Xsample{\bfX_s}

\def\Rsyn{{R_{\rm syn}}\xspace}
\def\musynr{{\mu_r^{\rm syn}}\xspace}
\def\sigsynr{{\sigma_r^{\rm syn}}\xspace}
\def\xisynr{{\xi_r^{\rm syn}}\xspace}

\def\Rerr{{R_{\rm bae}}\xspace}
\def\muerrell{{\mu_\ell^{\rm bae}}\xspace}
\def\sigerrell{{\sigma_\ell^{\rm bae}}\xspace}
\def\xiErrEll{{\xi_\ell^{\rm bae}}\xspace}

\definecolor{matlab1}{rgb}{0.0000, 0.4470, 0.7410}
\definecolor{matlab2}{rgb}{0.8500, 0.3250, 0.0980}
\definecolor{matlab3}{rgb}{0.9290, 0.6940, 0.1250}
\definecolor{matlab4}{rgb}{0.4940, 0.1840, 0.5560}
\definecolor{matlab5}{rgb}{0.4660, 0.6740, 0.1880}
\definecolor{matlab6}{rgb}{0.3010, 0.7450, 0.9330}
\definecolor{matlab7}{rgb}{0.6350, 0.0780, 0.1840}

\definecolor{matplotlib0}{HTML}{1f77b4}
\definecolor{matplotlib1}{HTML}{ff7f0e}
\definecolor{matplotlib2}{HTML}{2ca02c}
\definecolor{matplotlib3}{HTML}{d62728}
\definecolor{matplotlib4}{HTML}{9467bd}
\definecolor{matplotlib5}{HTML}{8c564b}
\definecolor{matplotlib6}{HTML}{e377c2}
\definecolor{matplotlib7}{HTML}{7f7f7f}
\definecolor{matplotlib8}{HTML}{bcbd22}
\definecolor{matplotlib9}{HTML}{17becf}

\ifpdf
\hypersetup{
  pdftitle={Multiway Synergy Emulator},
  pdfauthor={E. Newman, A. Arnold, D. Calvetti, and E. Somersalo}
}
\fi

\begin{document}

\maketitle

\begin{abstract}
The muscle recruitment problem in biomechanics characterizes all possible muscle activation patterns that correspond to an observed locomotion.  Due to the inherent redundancy of the musculoskeletal system, the ensemble of possible activations defines a high-dimensional manifold that can be explored by stochastic sampling methods such as Markov chain Monte Carlo (MCMC). It has been hypothesized that the central nervous system (CNS) does not control each muscle separately, but for standard movements such as level walking, key muscles are activated in few synergetic groups, using the redundancy for fine tuning the motion and balance. 
In this work, we propose a low-rank tensor decomposition approach to extract muscle activation pattens from MCMC sample data that corresponds to given observed dynamics.  The computed low-rank features capture the dominant dynamics and exhibit lower variability that the full MCMC realization of a specific component.  This corroborates the current uncontrolled manifold hypothesis of how the CNS controls the skeletomuscular system. We further leverage the tensor approximations to design an efficient low-rank emulator of the MCMC sampler. We demonstrate the approximation quality and computational efficiency of the tensor synergy extraction and emulation strategies using statistical, storage, and timing metrics. 
\end{abstract}



\section{Introduction}
\label{sec:introduction}

Human and vertebrate animal locomotion requires a complex concerted control of the numerous muscle forces over the time span of the motion tasks. Understanding the neuronal muscle activation patterns during even relatively simple tasks like plain walking and balance control remains one of the central questions in neuroscience and biomechanics. The skeletal muscles responsible for motion and balance of humans and vertebrate animals are an abundantly redundant system, with several significantly different muscle activation patterns resulting in essentially the same motor outcome. How, exactly, the central nervous system (CNS) selects the neural activation pattern of muscles, and how the feedback system works to regain balance and movement during a motor task subject to external perturbations is not fully understood. One popular hypothesis is that to generate the overall motion pattern, the CNS controls {\em synergetic muscle groups} rather than individual muscles,  using the extra degrees of freedom for fine tuning the balance and performance. An alternative hypothesis assumes that the muscle control is an expression of optimal control with a minimum intervention principle. This hypothesis is closely related to the uncontrolled manifold hypothesis \cite{scholz_uncontrolled_1999} 
and to the degree of freedom problem introduced and discussed earlier by Bernstein
\cite{bernstein_coordination_1967}. In this work, we refer to muscle synergies as a concept of concerted action of the muscles during the task, without adhering to any particular assumption of the role of the CNS.

The identification of putative muscle synergies is usually based on low-rank approximations through
different factorization methods of multi-muscle electromyography (EMG) time series data recorded during motor tasks \cite{safavynia_muscle_2011}. Two of the main drawbacks of the EMG data are that the number of muscles monitored during the task is limited, and that the voltage signal is very noisy and difficult to translate into muscle force scale. To circumvent these limitations, a potentially viable, non-invasive  alternative to EMG-based identification of the synergies was recently proposed in \cite{calvetti_separable_2025}, where the starting point for synergy analysis was the solution of {\em muscle recruitment problem} based on kinematic data in the Bayesian framework. The aim of the muscle recruitment problem is to identify  muscle activation patterns that effectuate an observed kinematic motion within the uncertainty of the data. The problem is ill-posed in the sense that there are infinitely many possible activation patterns that can be associated with the observed movement data; i.e., the number of essential degrees of freedoms to describe the motion of the body is significantly lower than the number of muscles.   A natural way to reduce the uncertainty inherently associated with the muscle recruitment problem is to design algorithms leveraging prior conditions excluding or discouraging activation patterns known to be physiologically less viable.  This can be achieved within the Bayesian framework for inverse problems by introducing physiology-informed priors. It is also important that such priors respect the redundancy, which is a key feature of the complex adaptive feedback system for motion and balance correction under changing conditions. In line with these ideas, a Markov chain Monte Carlo (MCMC) algorithm to draw a large and representative sample of time series of muscle activation patterns from the posterior distribution integrating the kinetic data and the prior model through a computational model was proposed in  \cite{amankwah_exploring_2024,bersani_myobolica_2024}, extending earlier work \cite{martelli_computational_2013,martelli_stochastic_2015}. Recently \cite{calvetti_separable_2025}, a synergy analysis for plain walk was carried out employing a low rank approximation of the posterior mean patterns generated by the MCMC sampler. 

In this paper, the MCMC sample is viewed as an third-order multiway array or \emph{tensor}, capturing the high-dimensional correlations between muscles, activation pattern time series, and sample variability. A low-rank tensor approximation is employed to identify muscle synergies, in accordance with the understanding that only a few of them essentially determine the motor task. Unlike in \cite{calvetti_separable_2025}, where the activation pattern based synergy analysis relied on the posterior mean of the activation patterns, our proposed synergy analysis is based on the full MCMC sample which provides a much richer set of posterior realizations characterizing the full posterior distribution of the activation patterns.  Given the stochastic nature of the MCMC sample, it is natural to ask how the posterior uncertainties are distributed among the synergies, and furthermore, which portion of the posterior uncertainty can be attributed to muscle activations that are not explained through the synergies. The latter question is particularly relevant for understanding how much leeway the CNS has for feedback-based trajectory corrections beyond the kinematics controlled by the putative synergy patterns. In this article, this question is answered  by applying low-rank tensor-based decomposition techniques to the full MCMC sample, and by analyzing the approximation discrepancy according to the Bayesian approximation error (BAE) paradigm. Moreover, the tensor based interpretation of the sample leads to a novel and effective simulation protocol of the muscle activation patterns that emulate the full MCMC sampler.

\subsection{Our Contributions}
\label{sec:contributions}

In this paper, leveraging the power of tensor analysis,  we propose a new methodology for a sample-based analysis of the muscle synergy. More precisely, in \Cref{sec:nncp}, we interpret the MCMC muscle recruitment sample as a multiway tensor and we employ a low rank nonnegative tensor decomposition strategy to extract interpretable features corresponding to random realizations of the muscle synergies. This interpretation, which to the best of our knowledge has not been previously described in the literature, provides a natural way to separate the uncertainties associated with the synergies from those related to genuine redundancies in the muscle control. 

Another novelty of our contribution is presented in \Cref{sec:cp_bae}, where the tensor analysis is carried out one step further and utilized to analyze portion of the sample not accounted for by the synergies in the framework of Bayesian approximation error (BAE). A byproduct of this interpretation is a novel tensor emulator pipeline for generating independent samples of muscle recruitment patterns that replicate the posterior distribution. Numerical experiments  in \Cref{sec:numerical_experiments} illustrate how synergies can be extracted from the MCMC sample through the tensor analysis and analyze the quality of the emulator, in terms of emulation metrics, reduced storage and computational cost. 

\subsection{Outline}
\label{sec:outline}

The paper proceeds as follows. 
In \Cref{sec:background}, we describe the muscle recruitment problem as a dynamic Bayesian linear inverse problem (\Cref{sec:muscle_model}) and describe prior synergy extraction techniques through a nonnegative matrix factorization (\Cref{sec:myobolica_nmf}). 
In \Cref{sec:methodology}, we describe our new synergy extraction methodology via a multiway nonnegative decomposition (\Cref{sec:nncp}) and a multiway emulation strategy inspired by the Karhunen–Lo\`{e}ve (KL) expansion (\Cref{sec:tensor_emulator}).  
In \Cref{sec:cp_bae}, we capture the variability across the samples using a Bayesian approximation error approach with an efficient multiway representation (\Cref{sec:tensor_bae}). 
In \Cref{sec:numerical_experiments}, we demonstrate the empirical performance of our proposed Mybolica emulator in terms of the quality of the extracted synergies (\Cref{sec:numerical_synergies}), the accuracy and efficiency of the full synergy and Bayesian approximation error emulation (\Cref{sec:numerical_emulation}). 
In \Cref{sec:conclusion}, we conclude and outline potential future work. 


\section{Background}
\label{sec:background}

In this section, we summarize the muscle recruitment problem, and its interpretation in the Bayesian framework. In this context the likelihood model  arises from the biomechanical model discussed in further detail in \Cref{sec:muscle_model_physics} of the appendix. Similarly, a brief review of the construction of the physiologically meaningful priors are given, postponing again the more technical  details to \Cref{sec:myobolica_prior} of the appendix.

The likelihood and the prior are combined according to Bayes' rule into posterior density. Given a set of kinetic data corresponding to a dynamic task, samples from its posterior can be generated efficiently by means of an MCMC sampler, Myobolica, specifically designed for these type of problems. For details of the sampler, we refer to \cite{bersani_myobolica_2024}. The posterior mean from the Myobolica sample was used to extract the synergies analyzed in 
\cite{calvetti_separable_2025}.

In the remainder of the paper, we will use lowercase letters for scalars, $x \in \Rbb$,  bold lowercase for vectors, $\bfx\in \Rbb^n$, bold uppercase for matrices, $\bfX\in \Rbb^{n_1\times n_2}$, and bold calligraphic symbols for tensors, $\TQ\in\Rbb^{n_1\times n_2\times n_3}$.  We reserve $t$ for time, $n$ with an informative subscript for particular dimensions, and $\pi$ with an informative subscript for either a probability density or distribution.

\subsection{Muscle Recruitment Model: A Dynamic Bayesian Linear Inverse Problem}
\label{sec:muscle_model}

For the sake of definiteness, we restrict the description here to the special case in which the dynamic task comprises few steps of level walk. The dynamics of the human walking is monitored by  following the trajectory of a set of fiducials attached to the body. The goal of the {\em inverse kinematic problem} in biomechanics is to solve for the time courses of the joint angles and the joint positions from the observed trajectories of the fiducials.
The related {\em inverse dynamics problem} seeks to estimate the associated joint torques from the joint angles and measured ground reaction forces.  To solve the inverse dynamics problem, one can model the body by a joint-link system and write the equations of motion in the Lagrangian dynamics. Effective tools for solving the  inverse kinematic and dynamics problems are available in computational frameworks such as OpenSim \cite{delp_opensim_2007}. As a starting point, we assume that the inverse kinematics problem and inverse dynamics problems have been solved to satisfactory accuracy. 

As explained in more detail in \Cref{sec:muscle_model_physics} of the appendix, the muscle forces and the joint torques are related to each other through a linear relation,
    \begin{align*}
     \bfb_0(t) = \bfA(t)\bfF(t), \quad 0\leq t\leq T,
    \end{align*}
where $\bfb_0\in\Rbb^d$ denotes the time dependent vector of the noiseless torques associated to the degrees of freedom at the joints capturing the geometry of the body, $\bfF(t) \in\Rbb^{n_m}$ is the vector of muscle forces associated to the $n_m$ muscles included in the kinetic model, and $\bfA(t) \in\Rbb^{d\times n_m}$ is the lever-arm matrix. Acknowledging that each muscle has an upper bound for the force it may generate, we scale the force components $F_j(t)$ by their maximum (tetanic) forces $F_j^*$, defining muscle activations $q_j(t)$ through $F_j(t) = F^*_j q_j(t)$. For simplicity, we assume that the maximum forces are independent of time. The muscle recruitment problem seeks to characterize the possible muscle activations $q_j(t)$ from the noisy reconstructions of the torques at the discrete times $t_k$, $1\leq k\leq n_t$. In matrix form, we write the model as
	\begin{align}
	\bfb(t_k) = \bfA(t_k)\bfF^\star \bfq(t_k) + \bfe(t_k)=\bfA^\star(t_k) \bfq(t_k) + \bfe(t_k), \quad 1\leq k\leq n_t
	\end{align}
where $\bfF^\star = {\rm diag}(F^\star_j)$, $\bfA^*$ is the scaled lever arm matrix, $\bfq(t_k) \in \Rbb^{\nmuscle}$ contains the muscle activations for each $\nmuscle$ muscle,
$\bfb(t_k)$ is the noisy reconstruction of the joint torques, and $\bfe(t_k)$, referred to as noise in the rest of the paper, accounts for the error in the computation of the torques. 
For additional details on the underlying biomechanics, see \Cref{sec:muscle_model_physics}. 

To set up the Bayesian model, we model $\bfe$  as scaled Gaussian noise,
$\bfe(t_k)\sim \Ncal({\bf0}_d, \sigma^2\bfI_d)$ with variance $\sigma^2 > 0$, yielding a likelihood model of the form
	\begin{align}
		\pi_{\rm lkh}(\bfb(t_k) \mid \bfq(t_k)) \propto \exp\left(-\tfrac{1}{2\sigma^2} \|\bfb(t_k) - \bfA^\star(t_k)  \bfq(t_k)\|_2^2\right).
	\end{align}
For later reference, we arrange the data in a matrix of size $m\times n_t$, denoting
  	\begin{align}
		\bfB = \begin{bmatrix}
			\vert & \vert & & \vert\\
			\bfb (t_1) & \bfb (t_2) & \cdots & \bfb (t_{\ntime})\\
			\vert & \vert & & \vert\\
		\end{bmatrix}. 
	\end{align}

Next we define a prior model capturing the physiological information concerning the muscle activations. The prior is based on three physiologically justifiable assumptions:
    \begin{itemize}
        \item[(a)] The muscle activations must respect the bound constraints $0\leq q_j(t)\leq 1$.
        \item[(b)] The CNS-muscle complex favors low activation patterns due to the minimum energy principle.
        \item[(c)] Under normal conditions, muscle activations avoid sudden muscle contraction patterns, i.e., patterns with low yank $\frac{d\bfF}{dt}(t) =\bfF(t)$ are favored \cite{lin_yank_2019}.
    \end{itemize}
To make the problem computationally feasible, the activation patterns are discretized in the time direction, and we define the $n_m\times n_t$-dimensional random variable,
	\begin{align}\label{Q matrix}
		\bfQ = \begin{bmatrix}
			\vert & \vert & & \vert\\
			\bfq (t_1) & \bfq (t_2) & \cdots & \bfq (t_{\ntime})\\
			\vert & \vert & & \vert\\
		\end{bmatrix}. 
	\end{align}
With this notation, the likelihood over all observation times is defined as
    \begin{align*}
     \pi_{\rm lkh}(\bfB\mid\bfQ) = \prod_{k=1}^d \pi_{\rm lkh}(\bfb(t_k)\mid \bfq(t_k)).
    \end{align*}
Denoting by $\pi_{\rm pr}(\bfQ)$ the prior density for $\bfQ$, it follows from Bayes' formula that the posterior density is
    \begin{align*}
     \pi(\bfQ\mid\bfB)\propto \pi_{\rm pr}(\bfQ)\pi_{\rm lkh}(\bfB\mid\bfQ).
    \end{align*}

Following \cite{martelli_computational_2013,martelli_stochastic_2015,amankwah_exploring_2024,bersani_myobolica_2024}, we explore the solution of the muscle recruitment problem by generating random samples from the posterior density. The computationally efficient Gibbs-within-Gibbs type  MCMC sampling algorithm {\em Myobolica}, introduced in \cite{amankwah_exploring_2024,bersani_myobolica_2024}, was designed specifically for the muscle recruiting problem. 

The starting point of the current analysis is a {\em Myobolica}-generated sequence of posterior realizations of muscle activations over time, $\{\bfQ^{(k)}\}_{k=1}^\nsample$. Each element of the sample is a matrix $\bfQ^{(k)}\in \Rbb_+^{\nmuscle \times \ntime}$ of the form 
	\begin{align}\label{eq:q_mat_sample}
		\bfQ^{(k)} = \begin{bmatrix}
			\vert & \vert & & \vert\\
			\bfq^{(k)}(t_1) & \bfq^{(k)}(t_2) & \cdots & \bfq^{(k)}(t_{\ntime})\\
			\vert & \vert & & \vert\\
		\end{bmatrix}. 
	\end{align}
The likelihood portion of the posterior accounts for correlations between different muscles, while the temporal correlations along the activation path of each muscle comes from the assumption (c) in the prior design. The spatio-temporal structure of the posterior sample makes it particularly malleable for multiway analysis, which is the core contribution of this article. Before presenting the details of that analysis, however, we give a brief summary of previous work on muscle synergies.

\subsection{Matrix-Based Myobolica Synergy Analysis}
\label{sec:myobolica_nmf}

Exploring the posterior via a Myobolica sample clearly displays the variability of muscle activation behavior during the motion task. Intuitively, the putative synergetic muscle control manifests itself in the form of inter-muscular and temporal correlations, making it possible to identify groups of muscles that follow similar temporal activation patterns.  A promising strategy introduced in \cite{calvetti_separable_2025}, potentially competitive over the traditional EMG analysis, extracts the muscle synergies from the posterior mean of the Myobolica-generated samples using a nonnegative matrix factorization (NMF) \cite{gillis_nonnegative_2020}.  The specific pipeline followed in \cite{calvetti_separable_2025} extracts $R$ muscle synergies via an NMF of the mean muscle-time activation matrix, $\overline{\bfQ} = \frac{1}{\nsample} \sum_{k=1}^\nsample \bfQ^{(k)}$; that is, 
	\begin{align}\label{eq:nmf}
		\overline{\bfQ} \approx \overline{\bfW}_m \, \overline{\bfW}_t^\top 
		= \sum_{r=1}^R \overline{\bfW}_m(:,r) \,\overline{\bfW}_t(:,r)^\top
	\end{align}
where $\overline{\bfW}_m\in \Rbb_+^{\nmuscle \times R}$ and $\overline{\bfW}_t\in \Rbb_+^{\ntime \times R}$ are the muscle and time factor matrix, respectively.  The $r$-th muscle synergy is the interpretable rank-one matrix $\overline{\bfW}_m(:,r) \,\overline{\bfW}_t(:,r)^\top$, where $\overline{\bfW}_m(:,r)$ contains the weights per muscle and $\overline{\bfW}_t(:,r)$ captures the activation pattern over time.  In the biomechanics literature, with reference to a similar decomposition of the EMG data, the former vectors are referred to as {\em synergy weights}, and the latter ones as {\em synergy excitation primitives}.
The product of the two encompasses the firing patterns of key muscle groups for the given task, e.g., human walking.  In sum, the decomposition of the mean activation into synergies is informative, and points towards the hypothesis that the muscle recruitment problem can be addressed by analyzing few muscle groups at a time.

\section{Multiway Analysis}
\label{sec:methodology}

While NMF has proven effective at synergy extraction, representing the posterior sample with its mean, as in \cite{calvetti_separable_2025},  effectively ignores  the variability over the samples, making it challenging to quantify the uncertainties of the synergy analysis. In the following, we propose an approach that respects the variability of the samples and exploit correlations along the entire stochastic process, using a nonnegative tensor decomposition. 

\begin{figure}
\centering
\scriptsize
 	\begin{subfigure}[b]{0.48\linewidth}
		\centering
		\includegraphics{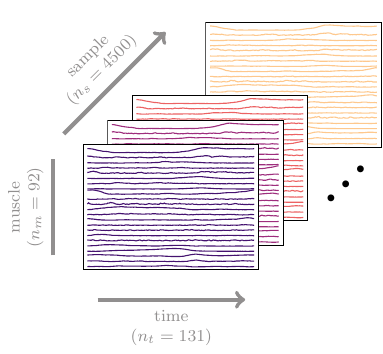}
		\caption{Example muscle-time activation patterns for various MCMC samples.  Each color corresponds to a different frontal slice/sample, $\TQ(:,:,k) = \TQ^{(k)}$.}
		\label{fig:data_visualization_3d}
	\end{subfigure}
	\hfill
	\begin{subfigure}[b]{0.48\linewidth}
		\centering
		\includegraphics{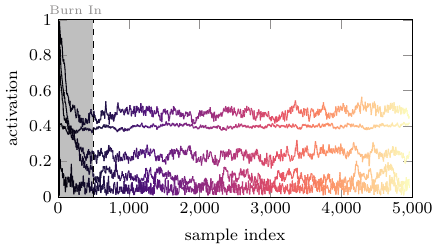}
		\caption{Example Markov chains for various muscle-time locations, e.g., $\operatorname{vec}(\TQ(i,j,:))\in \Rbb^{\nsample}$ for muscle $i$ and time $j$. Color corresponds to frontal slice/sample.}
		\label{fig:data_visualization_chains}
	\end{subfigure}
	
\caption{Visualization of Myobolica activation tensor, $\TQ \in \Rbb_+^{\nmuscle \times \ntime \times \nsample}$. }
\label{fig:data_visualization}
\end{figure}

To this end, we represent a sequence of $n_s$ muscle-time samples 
$\{\bfQ^{(k)}\}_{k=1}^\nsample$ generated by Myobolica as a multiway array or \emph{tensor}, $\TQ \in \Rbb_+^{\nmuscle \times \ntime \times \nsample}$, with dimensions corresponding to time, muscle, and sample, respectively.  Because of the three distinct dimensions, we call $\TQ$ a third-order tensor.  Using {\sc Matlab} indexing notation, the $k$-th frontal slice, $\TQ(:,:,k)\in \Rbb_+^{\nmuscle \times \ntime}$ is the $k$-th muscle-time matrix generated by the MCMC algorithm, i.e., $\TQ(:,:,k) = \bfQ^{(k)}$ for $k=1,\dots,\nsample$.  We visualize $\TQ$ in 
\Cref{fig:data_visualization_3d} and example Markov chains in \Cref{fig:data_visualization_chains}.

\subsection{Multiway Synergy Extraction}
\label{sec:nncp}

\begin{figure}
\centering
		\includegraphics{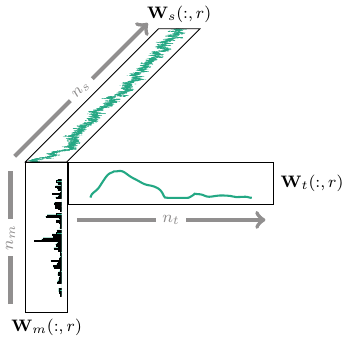}
        
		\caption{Rank-one synergy tensor $\Wmuscle(:,r) \circ \Wtime(:,r) \circ \Wsample(:,r)$. We depict the first mode, $\Wmuscle(:,r)$, identified as a synergy weight vector, using a bar chart. The largest values in magnitude indicate the muscles that are most active for this pattern. We depict the second mode, $\Wtime(:,r)$, or synergy activation primitive, as a continuous curve. The curve indicates the activation pattern over time, with larger values indicating that the muscles of the synergy group were activated most during a particular time period. The activity during the first half and inactivity during the second half is typical of the bipedal motion (i.e., left and right legs alternate activity). We depict the third mode, $\Wsample(:,r)$, using a trace plot. Each point represents the weight of the muscle-time pattern in a particular sample, and the stochasticity is expected due to the MCMC sampler. }
		\label{fig:synergy_rank_one_tensor}
\end{figure}

Motivated by NMF in \eqref{eq:nmf}, we leverage an analogous nonnegative tensor factorization based on the canonical polyadic/parallel factor (CP) decomposition \cite{hitchcock_expression_1927, carroll_analysis_1970, harshman_foundations_1970}.  We approximate the Myobolica tensor $\TQ$ by a rank-$\Rsyn$ nonnegative CP decomposition (NNCP), $\widehat{\TQ}_{\Rsyn} \in \Rbb_+^{\nmuscle \times \ntime \times \nsample}$, using a factor matrix per mode, concisely expressed in bracket notation \cite{kolda_tensor_2009} by 
    \begin{align}\label{eq:Q_approx}
		\TQ \approx \widehat{\TQ}_{\Rsyn} = \llbracket \Wmuscle, \Wtime, \Wsample \rrbracket \coloneqq \sum_{r=1}^\Rsyn \Wmuscle(:,r) \circ \Wtime(:,r) \circ \Wsample(:,r). 
	\end{align}
The muscle, time, and sample factor matrices are
$\Wmuscle \in \Rbb_+^{\nmuscle \times \Rsyn}$,  $\Wtime \in \Rbb_+^{\ntime \times \Rsyn}$, and  $\Wsample\in \Rbb_+^{\nsample\times \Rsyn}$, respectively.  Each summand is an  NNCP-rank-one muscle synergy tensor (\Cref{fig:synergy_rank_one_tensor})  and the symbol ``$\circ$" denotes a higher-order generalization of the outer product\footnote{The outer product of vectors $\bfa\in \Rbb^{n_1}$, $\bfb\in \Rbb^{n_2}$, and $\bfc\in \Rbb^{n_3}$  produces a tensor $\bfa\circ \bfb\circ \bfc \in \Rbb^{n_1\times n_2\times n_3}$ that is defined entry-wise by $(\bfa \circ \bfb \circ \bfc)(i,j,k) = \bfa(i) \bfb(j) \bfc(k)$ for each $i=1,\dots,n_1$, $j=1,\dots,n_2$, and $k=1,\dots,n_3$.}.
One can interpret \eqref{eq:Q_approx} as a simultaneous, weighted NMF decomposition of each frontal slice 
	\begin{align}\label{eq:weighted_nmf}
		\TQ(:,:,k) \approx \sum_{r=1}^\Rsyn \Wsample(k,r) \cdot (\Wmuscle(:,r) \Wtime(:,r)^\top) 
	\end{align}
for $k=1,\dots,\nsample$.  The muscle and time factors are shared among samples and $\Wsample(k,r) > 0$ is the weight given to $r$-th muscle-time synergy in the $k$-th sample. 

We construct the NNCP approximation by solving the constrained optimization problem 
    \begin{align}\label{eq:nncp_goal}
    \begin{split}
        &\min_{\Wmuscle, \Wtime, \Wsample} \left\|\TQ - \llbracket \Wmuscle, \Wtime, \Wsample \rrbracket \right\|_{F}^2\\
        &\text{s.t.} \qquad \Wmuscle \in \Rbb_+^{\nmuscle \times \Rsyn}, \Wtime \in \Rbb_+^{\ntime \times \Rsyn}, \Wsample \in \Rbb_+^{\nsample \times \Rsyn}
        \end{split}
    \end{align}
where $\|\cdot \|_{F}$ is the Frobenius norm\footnote{The Frobenius norm of a third-order tenosr is $\|\TX\|_{F}^2 = \sum_{i=1}^{n_3} \sum_{j=1}^{n_2} \sum_{k=1}^{n_3} |\TX(i,j,k)|^2$ for any $\TX\in \Rbb^{n_1\times n_2\times n_3}$.}. Algorithmically, we update each factor matrix cyclically using a multiplicative updating strategy from  \cite{lee_algorithms_2000} on various matricizations of the objective function. 
For example, to update the muscle factor matrix $\Wmuscle$, NNCP re-frames \eqref{eq:nncp_goal} as a nonnegative least squares problem
	\begin{align}\label{eq:nnls}
		\min_{\Wmuscle \in \Rbb_+^{\nmuscle\times \Rsyn}} \|\bfQ_{(1)} - \Wmuscle \bfH_m^\top \|_{F}
	\end{align}
where $\bfQ_{(1)} = \begin{bmatrix} \TQ(:,:,1)  & \cdots & \TQ(:,:,\nsample) \end{bmatrix} \in \Rbb_+^{\nmuscle \times (\ntime \nsample)}$ is the mode-$1$ unfolding of $\TQ$.  The helper matrix, $\bfH_m\in \Rbb_+^{(\ntime\nsample) \times \Rsyn}$, is built from the Khatri-Rao product, $\odot$, of the remaining factor matrices
	\begin{align}\label{eq:khatri_rao}
		\bfH_m 
			= \Wsample \odot \Wtime 
			\coloneq \begin{bmatrix}\vert &  & \vert\\ \Wsample(:,1) \otimes \Wtime(:,1) &  \cdots &  \Wsample(:,\Rsyn) \otimes \Wtime(:,\Rsyn)\\ 
            \vert &  & \vert    
            \end{bmatrix}
	\end{align}
where $\otimes$ is the Kronecker product \cite{petersen_matrix_2008}. One then applies a standard multiplicative update via
	\begin{align}\label{eq:multiplicative_update}
		\Wmuscle \gets \Wmuscle * \left(\dfrac{\bfQ_{(1)} \bfH_m}{ \Wmuscle\bfH_m^\top \bfH_m + \varepsilon \mathbf{1}_{\nmuscle \times (\ntime\nsample)}}\right)
	\end{align}
where $*$ denotes the Hadamard pointwise product, division is performed elementwise, $\mathbf{1}_{n_1\times n_2}$ is an $n_1\times n_2$ matrix of all ones,  and $\varepsilon > 0$ is a small additive perturbation of each entry to avoid division by zero. 
Similar updating strategies are then applied to $\Wtime$ and $\Wsample$ with different matricizations of $\TQ$. The process is repeated until convergence. We outline the NNCP algorithm in \Cref{alg:cp_nmu}.

\subsection{Tensor-based Separable Expansion}
\label{sec:tensor_emulator}

The extraction of the synergies through the multi-way decomposition discussed in the previous section reveals the low-dimensional hidden structure of the muscle activation sample. The synergies, however, do not fully account for the sample variability, which is important in understanding to what extent the activation patterns are task-specific, and the extent of the uncontrolled redundancy. We motivate the proposed tensor-based separable expansion from the Karhunen–Lo\`{e}ve (KL) expansion. The expansion yields naturally a low-cost tensor emulator of the MCMC sampler that could be a valuable tool for different downstream tasks in movement analysis.

We first frame the muscle activation patterns as an $\Rbb^{n_m}$-valued stochastic process. Given a probability space\footnote{A probability space is a triple consisting of a sample space $\Omega$, a Borel $\sigma$-algebra $\mathscr{B}$ of measurable sets, and probability measure $\Pbb$.} $(\Omega,{\mathscr B},\Pbb)$, let $\bfq(t; \cdot):\Omega \to \Rbb^{n_m}$ denote a square integrable stochastic process\footnote{Specifically, we have $\int_{\Omega} \int_0^\top \|\bfq(t; \omega)\|_2^2 \,\text{d}t\, \text{d}\Pbb(\omega) < \infty$.} over $t\in [0,T]$.  Square integrable random processes admit a Karhunen-Lo\'{e}ve expansion (see, e.g., \cite{sullivan_introduction_2015}) given by
    \begin{align}\label{eq:kl_expansion}
      \bfq(t; \omega) = \overline\bfq(t) + \sum_{j=1}^\infty  \xi_j(\omega)  \sqrt{\lambda_j} \bfu_j(t),
    \end{align}
where $\overline\bfq: [0,T] \to \Rbb^\nmuscle$ is the mean trajectory, $\lambda_j\geq 0$ are the eigenvalues of the covariance operator of $\bfq$, $\lambda_j \searrow 0$ as $j\to\infty$, and the functions $\bfu_j: [0,T] \to \Rbb^\nmuscle$ are mutually orthogonal, that is, $ \int_0^\top \bfu_j(t)^\top \bfu_k(t)\,\text{d}t = \delta_{jk}$ where $\delta_{jk}$ is the Kronecker delta. The real valued random variables, $\xi_j: \Omega \to \Rbb$, have zero mean, unit variance, and are mutually uncorrelated, that is,
    \begin{align}\label{eq:random_variable_uncorrelated}
      \Ebb\big[\xi_j\xi_k\big] = \int_\Omega \xi_j(\omega)\xi_k(\omega)\, \text{d} \Pbb(\omega) = \delta_{jk}.
    \end{align}
When available, the KL decomposition provides a tool for generating {\em independent} samples of $\bfq$ by generating univariate samples of the variables $\xi_j$. To simplify the draws further, it is common to approximate the variables $\xi_j$ by standard normal variables. In contrast to the multiway synergy analysis, however, the orthogonal functions $\bfu_j$ in the KL decomposition are usually not separable, i.e., they cannot be written as products of weight vectors and activation primitives, $\bfu_j(t) = \bfv_j g_j(t)$, where $\bfv_j\in\Rbb^\nmuscle$ is time independent.  In the following, we show that it is possible to find a separable approximation of the activation patterns based on tensor factorizations that allow independent sampling akin to the KL expansion. The derivation is done in two steps, first approximating the process by a low-rank NNCP decomposition, and subsequently using the Bayesian approximation error (BAE) analysis to process the approximation error.

\begin{figure}
\centering
		\includegraphics{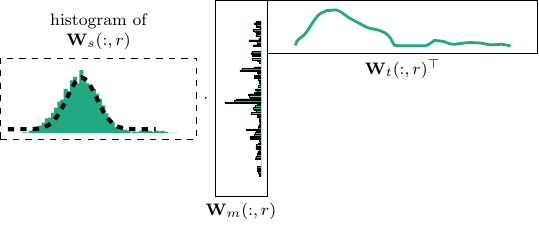}
		
		\caption{Rank-one synergy to emulator.  The $r$-th synergy, $\Wmuscle(:,r) \Wtime(:,r)^\top$, is retained from the NNCP. The black dashed curve is the Gaussian approximation formed from the empirical mean, $\musynr$, and variance, $(\sigsynr)^2$,  of the sample weights, $\Wsample(:,r)$.}
		\label{fig:synergy_rank_one_tensor_emulation}
\end{figure}

We start by interpreting a Myobolica muscle-time sample path as discrete measurements of particular realization of the stochastic process; i.e., for $\omega_k\in \Omega$, one can express  $\bfQ^{(k)}$ in \eqref{eq:q_mat_sample} as 
	\begin{align}
	\bfQ(\omega_k) = \begin{bmatrix}
			\vert & \vert & & \vert\\
			\bfq(t_1; \omega_k) & \bfq(t_2; \omega_k) & \cdots & \bfq(t_{\ntime}; \omega_k)\\
			\vert & \vert & & \vert\\
		\end{bmatrix}.
	\end{align}
Following the frontal slice presentation in \eqref{eq:weighted_nmf},  the NNCP approximation of the Myobolica tensor allows one to write the $k$-th  sample as a weighted sum of rank-one matrices 
	\begin{align}
		\bfQ(\omega_k)  
			\approx \sum_{r=1}^\Rsyn \Wsample(k,r) \cdot (\Wmuscle(:,r) \Wtime(:,r)^\top), 
	\end{align}
reminiscent of the  KL expansion in \eqref{eq:kl_expansion} with the rank-one matrices $\Wmuscle(:,r) \Wtime(:,r)^\top$ mirroring the role of the basis functions $\bfu_j$, however, unlike in the KL expansion, the random part $\Wsample$ is not centered. The main idea here is to retain the physical properties of the ``bases'' $\Wmuscle$ and $\Wtime$ from the muscle activation synergies,  replacing the sample weights $\Wsample$ with appropriate random variables.  To this end, we consider the sample weights as realizations of underlying random variables $\{\xi_r^{\rm syn}\}_{r=1}^\Rsyn$, with probability distributions that approximately coincide with those of the underlying sample, see Figure \ref{fig:synergy_rank_one_tensor_emulation} for a schematic outline.  Thus, a tensor-emulated sample can be generated by setting
    \begin{align}\label{eq:tensor_emulator}
		\widetilde{\bfQ}(\omega) = \sum_{r=1}^\Rsyn \xi_r^{\rm syn}(\omega) \cdot (\Wmuscle(:,r) \Wtime(:,r)^\top). 
	\end{align}
For computational convenience, the random variables $\xi_r^{\rm syn}$ are chosen to be independent univariate Gaussian, i.e., $\xisynr\sim \Ncal(\musynr, (\sigsynr)^2)$, which preserve, up to a scalar multiple, the uncorrelated random variable property  \eqref{eq:random_variable_uncorrelated} in the KL expansion.  
The mean and the variance of  $\xi_r^{\rm syn}$ are
estimated from the sample factor matrix $\bfW_s \in \Rbb_+^{\nsample\times \Rsyn}$, setting	
    \begin{align}\label{eq:guassian_statistics_estimation}
		\musynr = \frac{1}{\nsample} \sum_{k=1}^\nsample \bfW_s(k,r) \qquad \text{and} \qquad 
		(\sigsynr)^2 = \frac{1}{\nsample - 1} \sum_{k=1}^\nsample (\bfW_s(k,r) - \musynr)^2,
	\end{align}
thus ignoring the possible cross correlations between the synergies.

We summarize the steps for the construction of the Myobolica tensor emulator and for the emulation procedure in the following algorithms, namely the training algorithm (\Cref{alg:tensor_emulator_training}) that extracts synergies and Gaussian distribution parameters and the deployment algorithm (\Cref{alg:tensor_emulator_deployment}) that emulates the synergetic part Myobolica samples.

\begin{algorithm}
\caption{Synergy Tensor Emulator Training}\label{alg:tensor_emulator_training}
    \begin{algorithmic}[1]
        \State {\bf Inputs:} Myobolica tensor $\TQ\in \Rbb_+^{\nmuscle \times \ntime\times \nsample}$, user-prescribed ranks $\Rsyn$
        \State {\bf Return:}  Emulator parameters $\big\{\Wmuscle \in \Rbb_+^{\nmuscle \times \Rsyn}, \Wtime \in \Rbb_+^{\ntime \times \Rsyn}, \{(\musynr, \sigsynr)\}_{r=1}^\Rsyn\big\}$
        \Statex
       	\State Compute rank-$\Rsyn$ NNCP approximation of $\TQ \approx \llbracket \Wmuscle, \Wtime, \Wsample \rrbracket$ \Comment{\cref{alg:cp_nmu}}
	\State Compute empirical statistics from columns  $r=1,\dots, \Rsyn$ of $\Wsample$ via \Comment{\cref{eq:guassian_statistics_estimation}}
		\begin{align*}
		\musynr = \frac{1}{\nsample} \sum_{k=1}^\nsample \Wsample(k,r)   \quad \text{and} \quad
		 (\sigsynr)^2 = \frac{1}{\nsample - 1} \sum_{k=1}^\nsample (\Wsample(k,r) - \musynr)^2
		\end{align*}
    \end{algorithmic}
\end{algorithm}

\begin{algorithm}
\caption{Synergy Tensor Emulator Deployment}\label{alg:tensor_emulator_deployment}
\begin{algorithmic}[1]
	\State {\bf Inputs:} Emulator parameters $\big\{\Wmuscle, \Wtime, \{(\musynr, \sigsynr)\}_{r=1}^\Rsyn\big\}$, number of samples $N$
	\State {\bf Return:} Emulated sample $\widetilde{\bfQ} \in \Rbb_+^{\nmuscle\times \ntime}$
		\Statex
		\State Initialize $\widetilde{\bfQ} = {\bf0}_{\nmuscle \times \ntime}$
		\For{$r=1,\dots, \Rsyn$} \Comment{parallelizable}
			\State Sample weight for synergy $r$ via $\xisynr \sim \Ncal(\musynr, (\sigsynr)^2)$
			\State Update $\widetilde{\bfQ} \gets \widetilde{\bfQ} + \xisynr \cdot \Wmuscle(:,r) \Wtime(:,r)^\top$
		\EndFor
\end{algorithmic}

\end{algorithm}

\section{Memory-Efficient Bayesian Approximation Error Emulation}
\label{sec:emulation}
\label{sec:cp_bae}

While the NNCP extracted synergies in \Cref{sec:nncp} account for the dominant muscle activation, the low-rank approximation does not fully capture the variability across the sample. This is in line with the physiological interpretation of the synergies.   While the muscle synergies are expected to capture the essential part of locomotion under normal conditions, the muscle system should allow for a wider range of less controlled alternatives. This variability may be activated in response to unexpected events, e.g., in case of balance loss due to changes in the normal ambient, requiring muscle contraction patterns that are not regulated by the synergies, including, e.g., cocontractions that do not alter the motion pattern but secure the joints.  The Bayesian approximation error approach, introduced in \cite{kaipio_statistical_2005,kaipio_statistical_2007} provides a flexible framework for investigating the discrepancy between the full sample dynamics and its low-rank approximation. Below, we give a concise summary of the basic ideas.

Consider a random variable ${X}$ that is approximated by a random variable $\widehat{ X}$. The discrepancy, or approximation error, is denoted by ${E}$,
    \begin{align*}
     {X} = \widehat{X} + ({X} - \widehat{X}) =
     \widehat{X} +{E}, 
    \end{align*}
Assuming that we can generate a sample $\{(\bfx^{(k)},\widehat\bfx^{(k)})\}_{k=1}^K$ from the joint probability distribution of the pair $({X},\widehat{X)}$, we form a sample $\{\bfe^{(k)}\}_{k=1}^K$ of the approximation error ${E}$ by defining  $\bfe^{(k)} = \bfx^{(k)} -\hat{\bfx}^{(k)}$. 
A popular, computationally efficient approximation, referred to as the {\em enhanced error model}, replaces ${E}$ by a Gaussian random variable, ${E} \approx \widehat{E}\sim{\mathcal N}(\overline\bfe,\bfC)$, where the mean and covariance are sample-based empirical estimates, $\bar{\bfe} = \frac{1}{K} \sum_{k=1}^K \bfe^{(k)}$ and  $\bfC = \frac{1}{K-1}\sum_{k=1}^K (\bfe^{(k)} - \bar{\bfe}) (\bfe^{(k)} - \bar{\bfe})^\top$.

To bring the ideas of the BAE in the present context, we express the Myobolica activation tensor as the sum of the NNCP approximation~\eqref{eq:Q_approx} and the leftover discrepancy,
	\begin{align}\label{eq:tensor_approx_err}
		\TQ = \widehat{\TQ}_{\Rsyn} +  \TE
	\end{align}
where $\TE\coloneqq \TQ -  \widehat{\TQ}_{\Rsyn} \in \Rbb^{\nmuscle \times \ntime \times \nsample}$ is the model error\footnote{The model error depends implicitly on the rank of the synergy approximation, $\Rsyn$. We omit the explicit dependence for notational convenience.}. In the BAE framework, we interpret $\TQ$ as samples of the  random variable of interest and $\widehat{\TQ}_{\Rsyn}$ as the corresponding samples of the approximate random variable.  In the following, after reviewing the BAE approximation in the traditional vectorized form, we propose an alternative that leverages the tensor formalism, leading to a computationally more efficient formulation.

\subsection{Matricized Bayesian Approximation Error}
\label{sec:matrix_bae}

We set up the stage for the Laplace approximation of the model error by first vectorize the muscle-time activation samples using a mode-$3$ unfolding, $\bfE_{(3)} \in \Rbb^{\nsample \times (\nmuscle \ntime)}$. Following the standard BAE procedure (\cite{kaipio_statistical_2005,kaipio_statistical_2007}), we then compute the vectorized sample mean $\overline{\bfe} = \frac{1}{\nsample} \sum_{k=1}^\nsample \bfE_{(3)}(k,:)^\top$ and the corresponding empirical covariance matrix $\bfC \in \Rbb^{(\nmuscle \ntime) \times (\nmuscle \ntime)}$.  In this matricized formulation, simulating the approximation error involves sampling from an $(\nmuscle \ntime)$-dimensional Gaussian distribution. Storing and generating samples using the full covariance matrix is memory-intensive, requiring storage costs of $\Ocal((\nmuscle \ntime)^2)$; in our case, $\nmuscle \ntime = \Ocal(10^5)$. In practice, we can select a rank parameter $\Rerr \ll \nmuscle\ntime$ and form a one-time, offline low-rank approximation of the covariance matrix via the rank-$\Rerr$ truncated eigenvalue value decomposition, $\bfC \approx \bfU_{\Rerr} \bfLambda_{\Rerr} \bfU_{\Rerr}^\top$, and sample from $\Ncal(\overline{\bfe}, \bfU_{\Rerr} \bfLambda_{\Rerr} \bfU_{\Rerr}^\top)$.  The low-rank approximation reduces the storage cost of the covariance to $\Ocal(\nmuscle \ntime \Rerr)$ and the computational cost of generating samples is similarly reduced. The matricized BAE emulator generates vectorized samples via
	\begin{align}\label{eq:bae_matrix_emulator}
		\bfe(\omega) = \overline{\bfe} + \bfU_{\Rerr} \bfLambda_{\Rerr}^{1/2} \bfeta
	\end{align}
where $\bfeta \sim \Ncal({\bf0}_{\Rerr}, \bfI_{\Rerr})$ is standard Gaussian noise. We then reverse the vectorization and reshape $\bfe(\omega) \in \Rbb^{\nmuscle\ntime}$ into an $\nmuscle \times \ntime$ matrix to emulate the muscle-time pattern. 

\subsection{Multiway BAE Laplace Approximation Emulator} 
\label{sec:tensor_bae}
To take advantage of the tensor-based error structure and further reduce storage costs, we propose a new CP-based emulator as a multiway alternative to the vectorized Laplace approximation. 
Analogous to the matricized approach in \Cref{sec:matrix_bae}, we start by computing the mean of the approximation error over the samples, $\overline{\bfE} = \frac{1}{\nsample} \sum_{k=1}^\nsample \Ecal(:,:,k)$. Next, we form the mean-centered error tensor $\TE^{\rm c}$ by subtracting the mean muscle-time pattern from each frontal slice of the error, $\TE^{\rm c}(:,:,k) = \TE(:,:,k) -  \overline{\bfE}$ for $k=1,\dots, \nsample$.  We then approximate the mean-centered error tensor by an unconstrained CP decomposition
	\begin{align}\label{eq:cp_bae}
		\TE^{\rm c} \approx \widehat{\TE}^{\rm c}_{\Rerr} = \llbracket \Xmuscle, \Xtime, \Xsample\rrbracket
	\end{align}
where $\Xmuscle\in \Rbb^{\nmuscle \times \Rerr}$ is the BAE muscle factor matrix, $\Xtime\in \Rbb^{\ntime \times \Rerr}$ is the BAE time factor matrix, and $\Xsample\in \Rbb^{\nsample \times \Rerr}$ is the BAE sample factor matrix. Observe that for the model error, positivity constraint does not apply. 

As in the NNCP strategy, forming a CP decomposition amounts to solving an unconstrained least squares problem
    \begin{align}\label{eq:cp_goal}
        \begin{split}
        &\min_{\Xmuscle, \Xtime, \Xsample} \left\|\TE^{\rm c} - \llbracket \Xmuscle, \Xtime, \Xsample \rrbracket \right\|_{F}^2\\
        &\text{s.t.} \qquad \Xmuscle \in \Rbb^{\nmuscle \times \Rerr}, \Xtime \in \Rbb_+^{\ntime \times \Rerr}, \Xsample \in \Rbb_+^{\nsample \times \Rerr}.
        \end{split}
    \end{align}
This can be done via an alternating least squares (CP-ALS) strategy where each factor matrix is updated by solving a least squares subproblem.  For example,  $\Xmuscle$ can be updated by solving
	\begin{align}\label{eq:least_squares}
		\min_{\Xmuscle \in \Rbb^{\nmuscle\times \Rerr}} \|\bfE_{(1)}^{\rm c} - \Xmuscle \bfY_m^\top \|_{F}
	\end{align}
where $\bfE_{(1)}^{\rm c}$ is the mode-$1$ unfolding of $\TE^{\rm c}$ and $\bfY_m = \Xsample \odot \Xtime$ where $\odot$ is the Khatri-Rao product defined in \eqref{eq:khatri_rao}. 
Any direct or iterative least squares solver can be used for \eqref{eq:least_squares}.  The CP-ALS algorithm cycles through subproblems for various unfoldings of $\TE^{\rm c}$ until convergence or stagnation. Additional details about the procedure can be found in \Cref{alg:cp_als}. 

Following the procedure outlined in \Cref{sec:tensor_emulator}, \Cref{alg:tensor_emulator_training}, and \Cref{alg:tensor_emulator_deployment}, we construct a separable emulator for the Bayesian approximation error from the CP factors
	\begin{align}\label{eq:bae_emulator}
		\widetilde{\bfE}(\omega) = \overline{\bfE} + \sum_{\ell=1}^\Rerr \xiErrEll(\omega) \cdot (\Xmuscle(:,\ell) \Xtime(:,\ell)^\top)
	\end{align}
where $\{\xiErrEll\}_{\ell=1}^\Rerr$ is a set of real-valued independent Gaussian random variables with mean and variance 
	\begin{align}
		\muerrell = \frac{1}{\nsample} \sum_{k=1}^\nsample \Xsample(k,r)   \quad \text{and} \quad
		 (\sigerrell)^2 = \frac{1}{\nsample - 1} \sum_{k=1}^\nsample (\Xsample(k,r) - \muerrell)^2.
	\end{align}
In summary, combining the synergy extraction and the Bayesian approximation error emulators yields a {\em fully separable Gaussian approximation}	of the form              
    \begin{align}\label{eq:tensor_emulator_with_bae}
		\widetilde{\bfQ}^{\rm full}(\omega) = \sum_{r=1}^\Rsyn \xi_r^{\rm syn}(\omega) \cdot (\Wmuscle(:,r) \Wtime(:,r)^\top)
			+ \overline{\bfE} +  \sum_{r=1}^\Rerr \xiErrEll(\omega) \cdot (\Xmuscle(:,\ell) \Xtime(:,\ell)^\top). 
	\end{align}
This formula, similar to the classical KL decomposition, allows independent sampling of the activation patterns by simple random draws from independent univariate Gaussians.

\section{Numerical Experiments}
\label{sec:numerical_experiments}

Our numerical tests use the fifth dataset of the Knee Grand Challenge described in \cite{fregly_grand_2012}: motion capture and ground reaction data of level walking of an 86 year old man with one instrumented knee prosthesis, antropomorphically scaled generic Gait2392 model \cite{delp_interactive_1990} is processed, and the corresponding joint torques and level arm matrices are estimated by tools provided by OpenSim. The model of the lower part of the body (lower torso, hip and legs) comprises 92 muscles, 43 in each leg, and 17 degrees of freedom, corresponding to ankles, subtalar joints, metatarsiphangeal joints, knees, hip, as well as joints between the pelvis and the torso. The muscle-joint system of the upper body is not included in the analysis. The lever arm matrices and the joint torques were used as the input for the Myobolica MCMC sampler, and a small sample $\{\bfQ^{(k)}\}_{k=1}^K$ of size $K=5000$ was generated. By visual inspection, the first $500$ sample matrices were considered as a burn-in and discarded, after which the sampler seems to be well-mixing. The samples were saved in the tensor $\TQ\in \Rbb_+^{\nmuscle \times \ntime\times \nsample}$ where $\nsample = 4500$. Code to reproduce the results is available upon request. 

\subsection{Synergy Extraction}
\label{sec:numerical_synergies}

\begin{figure}
\centering
\includegraphics{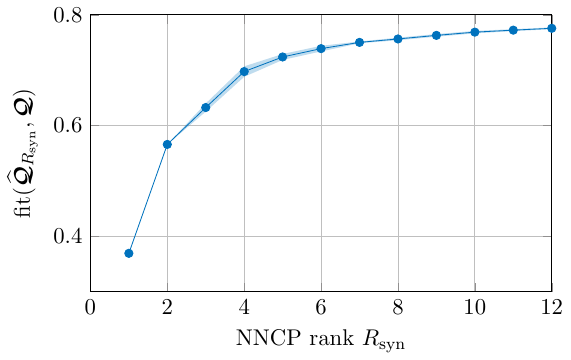}

\caption{Goodness of fit of the NNCP approximation $\widehat{\TQ}_\Rsyn$ to the Mybolica muscle activation tensor $\TQ$ for various ranks $\Rsyn = 1,\dots,12$.  The final fit score, $\fitmetric(\widehat{\TQ}_\Rsyn, \TQ) = 1 - \|\TQ - \widehat{\TQ}_\Rsyn\|_{F} / \|\TQ\|_{F}$, measures the fraction of the actual tensor captured by the low-rank NNCP approximation (higher is better).  Results were averaged over $10$ random algorithm initializations and results were reported using $4500$ samples from $\TQ$.  Nearly identical fit metrics  were obtained when using different $\TQ$ built from $1500$, $2500$, and $3500$ samples.}
\label{fig:mean_std_over_seeds}
\end{figure}

We use the {\sc Matlab} implementation in the Tensor Toolbox \cite{bader_matlab_2017} to compute the NNCP decomposition of the Mybolica muscle activation tensor. Following the Tensor Toolbox defaults, we randomly initialize the NNCP factors entrywise by drawing i.i.d. samples uniformly from the interval $[0.1, 1]$.  The NNCP implementation proceeds with the multiplicative updating strategy outlined in Algorithm~\ref{alg:cp_nmu} until a desired tolerance of $\delta = 10^{-4}$ is reached on the fit metric given by 
	\begin{align}
		\fitmetric(\widehat{\TQ}_\Rsyn, \TQ) \coloneqq 1 - \frac{\|\widehat{\TQ}_\Rsyn - \TQ\|_{F}}{\|\TQ\|_{F}}.
	\end{align}
The fit metric quantifies the fraction of $\TQ$ explained by the NNCP approximation, $\widehat{\TQ}_\Rsyn$. A larger score indicates a better fit, with a perfect fit corresponding to a value of one. We run our experiments various NNCP ranks and $10$ different random initializations per rank. The experimental results, displayed in Figure~\ref{fig:mean_std_over_seeds}, show the dependency of the fit on the rank parameter, $\Rsyn$.  As expected, the fit metric score increases with the NNCP rank, and the rate of increase slows at larger ranks, saturating right  below the value of $0.8$. This indicates that the underlying synergies can be well-approximated by a low-rank NNCP representation. The performance is consistent across different NNCP initializations. 

\begin{figure}
\centering
\includegraphics{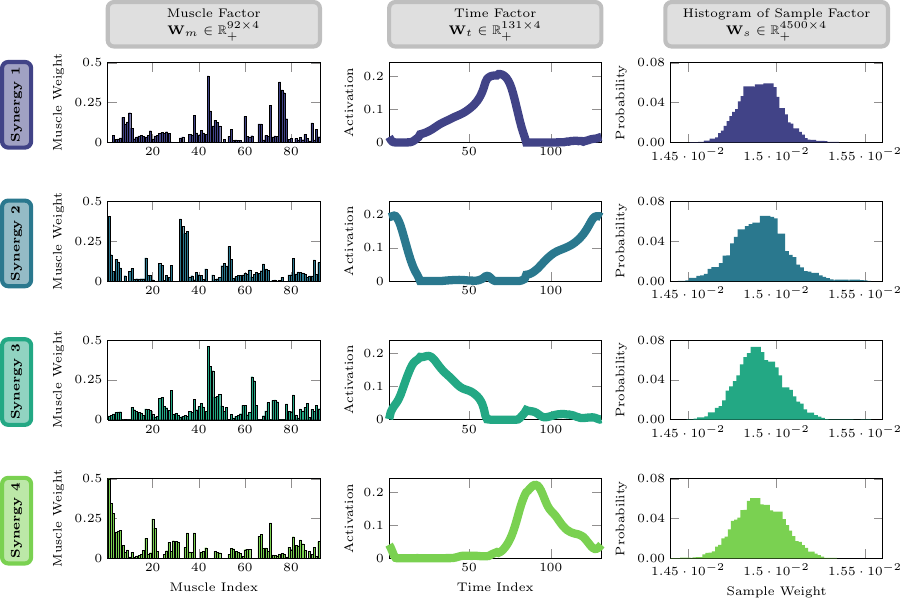}

\caption{Extracted NNCP synergies with CP rank $R = 4$ of the Myobolica muscle activation tensor $\TQ$. The panels in each column corresponds to a particular factor matrix: muscle in the first column, time in the second column, and sample in the third column. The muscle activation weights (first column) are displayed as bar plots,  each bar indicating to the weight given to the corresponding muscle. The panels in the second column display the activation primitives of the synergies over the time span of the motor task as continuous curves.  The third column displays the weights of the synergies across the sample in histogram format. The rows are organized according to the CP weights of the extracted synergies in decreasing order, i.e.,  $\lambda_1 \approx  449$, $\lambda_2 \approx 443$, $\lambda_3 \approx 380$, and $\lambda_4 \approx 371$, respectively. 
}
\label{fig:cp_factors}
\end{figure}

\begin{figure}
\centering

	\includegraphics{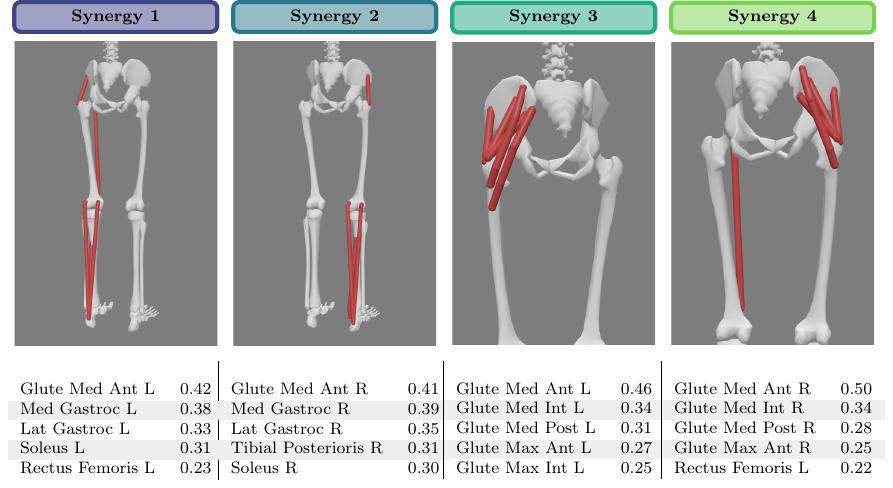}

    \caption{Muscle synergy plots of $\Wmuscle(:,r)$ for $r=1,\dots, \Rsyn$ generated in OpenSim. Below each image is a list of the top five weighted muscles per synergy. 
    We report the weight of each muscle to the right of the muscle name. Muscle names ending with ``L'' corresponds to muscles in the left leg and ending with ``R'' to the right leg. Abbreviations include Glute (Gluteus), Gastroc (Gastrocnemius), Lat (Lateral), Ant (Anterior), Int (Interior), Post (Posterior), and Max (Maximus).  The abbreviation Med stands for Medial if it appears before the muscle and Medius if it appears after the muscle.  Observe a little  break of the left-right symmetry, possibly due to the asymmetric gait caused by the prosthetic knee of the subject. 
    }
    \label{fig:muscle_synergy_skeleton}
\end{figure}

In Figure~\ref{fig:cp_factors}, we display the extracted synergies for NNCP rank $\Rsyn=4$. Figure~\ref{fig:muscle_synergy_skeleton} shows the anatomical location of the five dominant muscles in each synergy, listed below each image. We observe that most of the dominant muscles in Synergies $1$ and $3$ are in the left leg, while those in Synergies $2$ and $4$ are in the right leg, indicating a natural chiral parity of the synergies. Overall, the top weighted muscles in Synergies $1$ and $2$ are similar and in opposite legs, indicating a nearly symmetric gait simulated by the Myobolica sampler; this also holds for Synergies $3$ and $4$. The symmetry breaks are to be expected, as the subject has a knee prosthesis that is likely to cause an asymmetric gait. This left-right correspondence is reflected in the extracted patterns of the time-dependent activation primitives $\Wtime$ displayed in the middle column of Figure~\ref{fig:cp_factors}.  Here, we observe that when Synergy $1$ is active (high values), Synergy $2$ is inactive, and vice versa; the same pattern is observed between Synergies $3$ and $4$. The shape of the time courses, together with the function of the most significant muscles indicate that Synergies 1 and 3 can be associated with  the stance phase of the left leg, and Synergies 2 and 4 with the swing phase, while for the right leg, the roles are switched. Finally, the third column in Figure~\ref{fig:cp_factors} shows the histograms of weights in each synergy. Notably, the histograms are tightly concentrated around the mean value, with a variability only of the order of one percent. This finding suggests that the synergies do not account for most of the uncertainty observed in the Myobolica samples, confirming the intuition that most of the uncertainty in the muscle recruitment is uncontrolled by the synergistic activations.

\subsection{Emulation}
\label{sec:numerical_emulation}

We assess the performance full tensor emulator defined in \eqref{eq:tensor_emulator_with_bae} through three different studies, summarized in three figures.  In \Cref{fig:quantiles}, we illustrate how the tensor emulator can capture the shape and variability of the Myobolica activation training patterns through quantile visualization. In \Cref{fig:lda}, we show that the tensor emulator generates samples that well-approximate the true Myobolica samples, as assessed through linear discriminant analysis. In \Cref{fig:err_vs_storage}, we demonstrate the computational and storage efficiency of the tensor emulator. 

All three explorations compare tensor emulation with an NNCP synergy emulator \eqref{eq:tensor_emulator} and either a CP-based Bayesian approximation error emulator \eqref{eq:bae_emulator} or a matrix-based BAE emulator \eqref{eq:bae_matrix_emulator}. We generate $K = 5000$ emulated samples for various BAE ranks, $\Rerr$.   Let $\widetilde{\TQ} \in \Rbb^{\nmuscle \times \ntime \times 5000}$ denote the tensor of emulated samples. 

\subsubsection{Tensor Emulator Captures Physical Behaviors and Variability}

\begin{figure}
\centering

\includegraphics{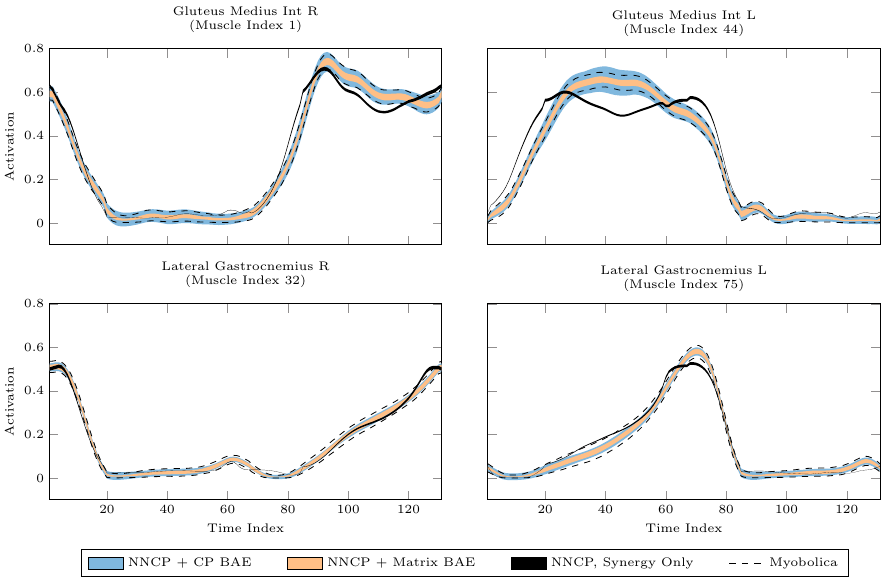}

\caption{Emulator approximations of the $90^{\rm th}$ quantile of four different muscles over time.  The first column shows the right muscle and the second column shows the left muscle. The first row shows the activation behavior of the Gluteus Medius Int and the second row shows the behavior for the Lateral Gastrocnemius. These were two highly weighted muscles in Synergies 1 and 2 in \Cref{fig:muscle_synergy_skeleton}. The black dashed lines show the $90^{\rm th}$ quantile from the Myobolica training data tensor, $\TQ$. Our goal is for the $90^{\rm th}$ quantiles of the emulators to match the black dashed boundaries. We use an NNCP emulator with $\Rsyn = 4$ in all emulations. The thin, solid black line shows the $90^{\rm th}$ quantile when we emulate with only the NNCP synergies from \eqref{eq:tensor_emulator}. While the synergy-only emulator captures the dominate behavior of each muscle, the low-rank estimate does not capture the variability of the Myobolica samples. The filled {\bf \color{matplotlib1}orange} region shows the  NNCP + Matrix BAE emulator using the matricized BAE from \eqref{eq:bae_matrix_emulator} with $\Rerr = 128$. The filled {\bf \color{matlab1}blue} region shows the $90^{\rm th}$ quantile obtained from NNCP + CP BAE emulator from \eqref{eq:tensor_emulator_with_bae} with $\Rerr = 128$. The full tensor emulator ({\bf \color{matlab1}blue}) more accurately captures the variability from the Myobolica $90^{\rm th}$ quantile than the matricized emulator ({\bf \color{matplotlib1}orange}) using more than 30 times less storage.  
}
\label{fig:quantiles}
\end{figure}

We first explore the ability of the tensor emulator to produce physically realistic estimates of the underlying Myobolica samples. In \Cref{fig:quantiles}, we begin by comparing the emulated samples to instances of the Myobolica sampler for a selected prominent muscle in the synergy extraction (\Cref{fig:muscle_synergy_skeleton}), the Gluteus Medius, the primary hip abductor and pelvic stabilizer during gait. We depict the $90^{\rm th}$ quantile of the emulated sample per muscle as follows. For the $i^{\rm th}$ emulated muscle, $\widetilde{\TQ}(i,:,:)\in \Rbb^{\ntime \times \nsample}$, we visualize the bounds of the interval containing $90^{\rm th}$ of the sample values per time point. We observe that the for the CP BAE emulator with $\Rerr = 128$ ({\bf \color{matlab1} blue}), the $90^{\rm th}$ quantile of specific muscle time courses well approximates the $90^{\rm th}$ quantile of the true Myobolica time courses (black dashed lines). Observe that the CP BAE version has larger quantile bounds than the matrix BAE version ({\color{matplotlib1} \bf orange}) with the same $\Rerr$. This indicates the all-CP emulator more effectively approximating the variability of the Myobolica samples. Note that if we  the $90^\text{\rm th}$ quantile the synergy emulator from the NNCP only, without the Bayesian approximation error ({\bf black}), we capture the dominant activation patterns, but not the variability. This further evidences that the synergies capture locomotion behavior under normal conditions, but are unable to capture  variability of gaits.

\subsubsection{Linear Discriminant Analysis} 
We compare the performance gains of the NNCP with the CP BAE emulator in \eqref{eq:tensor_emulator} over the  NNCP with the matrix BAE emulator \eqref{eq:bae_matrix_emulator}. We compare samples generated from each emulation strategy, $\widetilde{\TQ}$, to the Myobolica-generated samples, $\TQ$, using linear discriminant analysis (LDA),  a dimension reduction tool to assess how close two multivariate distributions are by comparing projected samples generated from each distribution.  In a nutshell, LDA determines a projection direction along which the within-class variance of the two samples to be compared is as small as possible (determine dominant  features of each set of samples), while the between-class variance is as large as possible (maximally separate samples from $\TQ$ and $\widetilde{\TQ}$). We perform LDA analysis per muscle, and describe the process in \Cref{sec:lda_details}. 

\begin{figure}
\centering

\includegraphics{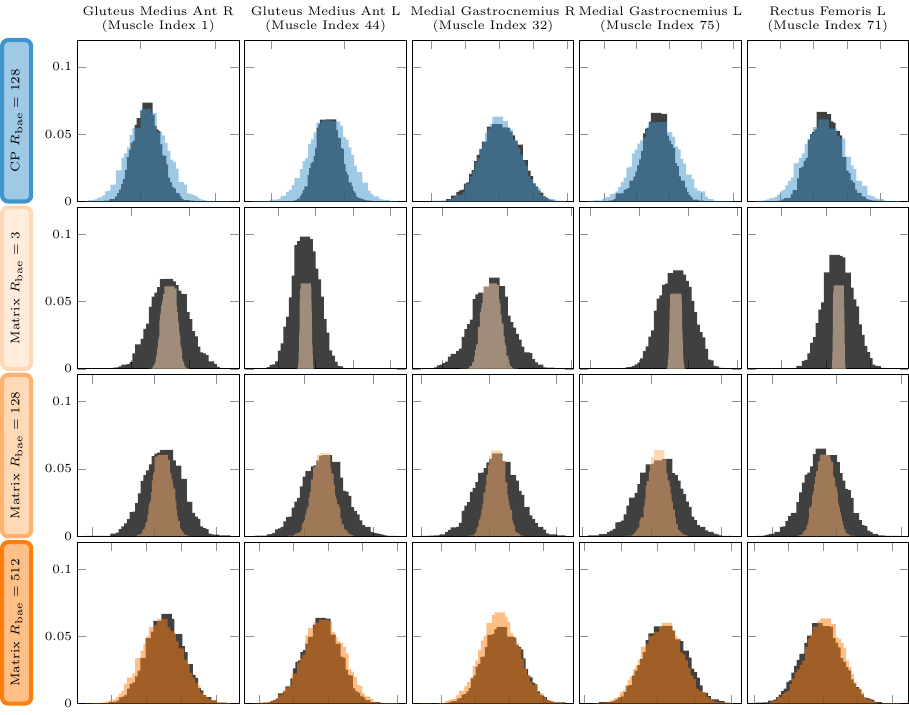}

\caption{Comparison of histograms via linear discriminant analysis (LDA) per muscle. The black histograms correspond to features obtained from the original Myobolica samples in $\TQ$. The colorful histograms correspond to features obtained from the $5000$ samples generated by the tensor emulator with either CP Bayesian approximation error (top row, {\bf \color{matlab1} blue}) with $R_{\rm cp} = 128$ or matrix approximation error (bottom three rows, {\bf \color{matplotlib1} orange}) for various ranks $R_{\rm mat} = 3$, $128$, and $500$. Observe that the CP emulator outperforms the matrix versions unless the rank is increased fourfold compared to that of the CP version.}

\label{fig:lda}
\end{figure}

To visually assess the quality of tensor factorization based emulator, we select six dominant muscles the synergies from \Cref{fig:muscle_synergy_skeleton}, and plot the histograms of the samples projected to the LDA direction, shown in columns of \Cref{fig:lda}. The underlying black histogram is the LDA projection of the Myobolica sample, and the colored histogram overlying it is the LDA projection of a selected emulated sample.

In the first row, the reference sample is compared to the emulator sample obtained by the CP tensor approximation \eqref{eq:bae_emulator}, with the CP tensor rank being $\Rerr = 128$. For comparison, the rows 2, 3 and 4 correspond to the vectorization-based matrix approximation \eqref{eq:bae_matrix_emulator}, computed by retaining $\Rerr = 3$, $\Rerr=128$, and $\Rerr=512$ eigenvectors of the large empirical covariance matrix.  We observe that the CP approximation (first row) outperforms the matrix approximations in quality, unless the number of retained eigenvectors is increased significantly (last row), however, the storage costs for retaining the eigenvectors in memory increases significantly. 

\begin{figure}
\centering

\includegraphics{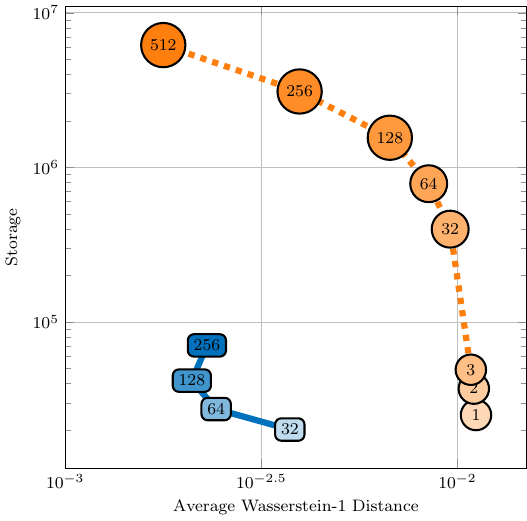}

\caption{Error vs. Storage (lower left is better). This diagram summarizes the trade-off between quality of the approximation by the Gaussian model \eqref{eq:bae_matrix_emulator} and by the CP tensor model \eqref{eq:bae_emulator}. The horizontal axis indicates the Wasserstein-1 distance between the LDA histograms, averaged over all muscles.   The vertical axis is the storage cost the particular emulator, computed by summing the relevant formulas in \eqref{eq:storage_computation}. The {\bf \color{matlab1}blue} curve indicates the performance of the CP approximation in this coordinate system, the numbers indicating the rank $R$ of the approximation, while the {\bf \color{matplotlib1}orange} curve refers to the performance of the Gaussian matrix approximation, the numbers indicating the rank of the covariance matrix approximation. We observe that the CP approximation emerges clearly superior in this accuracy-complexity comparison. }
\label{fig:err_vs_storage}
\end{figure}

\subsubsection{Accuracy-Complexity Tradeoff} 
In Figure~\ref{fig:err_vs_storage}, we compare the accuracy-complexity tradeoff by plotting the Wasserstein-1 distance of the LDA histograms, averaged over the muscles, versus the total storage of an emulator. The Wasserstein-1 distance approximates the absolute difference between empirical cumulative distribution functions. Storage is computed as the sum of the storage of the NNCP synergy emulator in \eqref{eq:tensor_emulator} and the storage of the BAE emulator, either for the CP BAE in \eqref{eq:bae_emulator} or the matrix BAE in \eqref{eq:bae_matrix_emulator}. The storage cost of each, denoted $\texttt{st}[\cdot]$,  is the total number of floating point numbers needed to form the various emulators; specifically,  
    \begin{subequations}\label{eq:storage_computation}
    \begin{align}
        \begin{split}
        \texttt{st}[\text{NNCP emulator}]
        &= \texttt{st}[\Wmuscle] + \texttt{st}[\Wtime] + \texttt{st}[\{\musynr, \sigsynr\}_{r=1}^\Rsyn]\\ 
        &=\nmuscle \Rsyn + \ntime \Rsyn + 2\Rsyn
        \end{split}\\
        \begin{split}
        \texttt{st}[\text{CP BAE emulator}]
        &=\texttt{st}[\Xmuscle] +\texttt{st}[\Xtime]  + \texttt{st}[\overline{\bfE}]+ \texttt{st}[\{\muerrell, \sigerrell\}_{\ell=1}^\Rerr]\\
        &=\nmuscle \Rerr+ \ntime \Rerr + \nmuscle\ntime  + 2\Rerr
        \end{split}\\
        \begin{split}
        \texttt{st}[\text{Matrix BAE emulator}]
        &= \texttt{st}[\Xmuscle] +\texttt{st}[\Xtime] + \texttt{st}[\overline{\bfE}] + \texttt{st}[\bfU_\Rerr] + \texttt{st}[\bfLambda_\Rerr] \\ 
        &=\nmuscle \Rerr+ \ntime \Rerr + \nmuscle\ntime + {\color{red}\nmuscle\ntime \Rerr} + \Rerr
        \end{split}
    \end{align}
    \end{subequations}
The key storage cost difference comes from the storage of the eigenvectors $\bfU_\Rerr$ in matrix BAE emulator, highlighted in red above. 

In \Cref{fig:err_vs_storage}, the results demonstrate that the emulator with the CP BAE outperforms the matrix BAE case in terms of accuracy-cost tradeoff. The matrix BAE case achieves a lower accuracy only for the highest rank, $\Rerr = 512$, and requires about 100 times more storage.  In contrast, the matrix BAE case requires less storage for the lower ranks, $\Rerr = 1, 2, 3$, but incurs about three times greater error.

\subsubsection{Computational Efficiency}
\label{sec:computational_efficiency}

\begin{figure}
\centering
\includegraphics{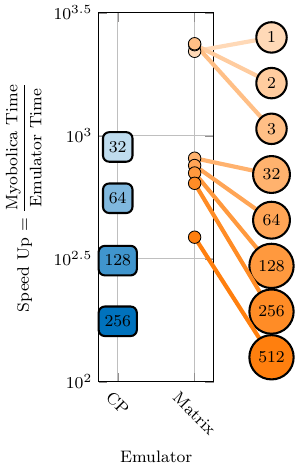}

\caption{Speed up, computed by the ratio of Myobolica sampler time to emulation time. We generate $50$ samples over $10$ independent runs.  In \Cref{fig:speed_up}, we report the average speedup over $10$ independent runs of generating $50$ samples with either the MCMC Myobolica sampler or the proposed emulation strategies.  Values greater than one indicate speedup and higher values are better.}
\label{fig:speed_up}
\end{figure}

As demonstrated by the discussion above, the tensor emulator produces physically realistic and statistically reliable estimates of the true Myobolica samples. In \Cref{fig:speed_up}, we show that the tensor emulator generates sample vectors with a speedup factor of 100-1000 over the original Myobolica sampler. Moreover, being an MCMC algorithm, the Myobolica sampler generates correlated samples, while the emulator samples are independent and therefore uncorrelated.  This observation is important if a small number of new sample vectors are desired, as the independently drawn samples reflect the full variability of the training set better. We note that the reported speedup does not include the training of the emulator, which we consider to be a one-time, offline cost.  We also note that the CP BAE emulator could be further accelerated by parallelizing sampling over the synergies.

\section{Conclusions}
\label{sec:conclusion}

This article proposes a novel way to extract muscle synergies from MCMC solution to the inverse muscle recruitment problem, exploring the muscle activation patterns based on the inverse kinematic solution of motion capture data. Unlike the previously proposed solution \cite{calvetti_separable_2025}
where the synergies were extracted from the posterior mean estimate only, here we exploit the tensor structure of the sampling data, and employ the low rank tensor approximation techniques to the data set. This way, we are able to rigorously assess the uncertainty of the synergies across the full sample without the need of projecting individual realizations on the synergy subspace. Furthermore, by combining the synergy analysis with a Bayesian approximation error (BAE) model, the tensor approximation provides a separable decomposition of the full sampling data, akin to the Karhunen-Lo\`{e}ve decomposition of multivariate stochastic processes, allowing us to emulate the muscle synergies by independent sampling. Combining the BAE further with low rank tensor decomposition, the full analysis leads to a memory efficient approximate stochastic model of the training sample, with the further advantage that it is possible to generate new independent samples from the underlying distribution without resorting to the MCMC scheme. While the analysis is this article is based on the muscle recruitment model, we point out that the same methodology allows the analysis and compression of Electromyography (EMG) data with essentially no necessary modification of the computational algorithms.

The uncertainty quantification analysis of this paper reveals that if the number of synergies in the level walking task is kept low, which is a reasonable assumption from the point of view of explaining the muscle control by the central nervous system under normal conditions, the uncertainties of the synergies explain only a fraction of the overall uncertainty identified by the MCMC algorithm \cite{amankwah_exploring_2024,bersani_myobolica_2024}. This finding seems to be in line with the intuition that uncertainties in the muscle activations serve as a rich resource of redundancy that allows the system for self-correcting when the synergetic control for some reason fails but the movement needs to be maintained.

\section*{Acknowledgments}
In addition to the funding acknowledgements on the first page, the authors acknowledge hardware and AI support. The numerical experiments were performed through the hardware made available in the Tufts University High Performance Compute Cluster \cite{phimmasen_reference_2021}. The authors also acknowledge the use of Claude Sonnet 5 in primarily for SLURM orchestration and to debug \LaTeX compilation errors, particularly to correct the \texttt{tikz} \cite{tantau_tikz_2023} externalize feature. Google Gemini AI Mode was used to suggest MSC codes. The authors take full responsibility for the accuracy and integrity of all content in this paper.

\appendix

\section{Underlying Biomechanics of Muscle Recruitment}
\label{sec:muscle_model_physics}
We briefly summarize the kinematic model underlying the analysis of the paper. Consider the human body, or a restricted part of it engaged in a physical task. We assume that it is characterized by $d$ degrees of freedom, typically joint angles, and we denote the time-dependent  joint torques by $\bfT_j(t)\in\Rbb^3$, $1\leq j \leq d$, where $0\leq t\leq T$. We assume that the body, or part of the body under consideration, is controlled by $n_m$ muscles, and we denote the muscle force vectors by $\bfF_k(t)\in\Rbb^3$, $1\leq k\leq n_m$. The torques at each joint can be calculated from the muscle forces, and since both the muscle force directions and the directions of the torque vectors are known quantities during the motion, the torque amplitudes can be related to the muscle force amplitudes through a known, time-dependent lever arm matrix. More precisely, let the torques and forces be written as $\bfT_j(t) = b_j(t) \widehat\bfT_j(t)$, $\bfF_j(t) = f_j(t) \widehat\bfF_j(t)$, where  $\|\widehat\bfT_j\| = \|\widehat\bfF_j\| = 1$, and since the geometry of the body is known, the direction vectors are known. On the other hand, the torque generated by a muscle at a given joint is related to the force through the lever arm relations $\bfT = \bfr\times\bfF$ at each joint, where $\bfr$ is the lever arm vector pointing from the pivot point of the joint to the attachment point of the muscle. The resultant torque at each joint is the sum of the torques of the attached muscles, and we can arrange these relations as a linear system relating the torque amplitudes and force amplitudes,
    \begin{align}
     \bfb(t) = \bfA(t) \bff(t) + \bfe(t), \quad 0\leq t\leq T,
    \end{align} 
where $\bfA(t)\in\Rbb^{d\times n_m}$ is the lever arm matrix depending on the geometry of the body, $\bfb(t) \in\Rbb^d$ and $\bff(t) \in\Rbb^{n_m}$ are vectors  with components $b_j(t)$ and $f_j(t)$, respectively, and $\bfe(t)$ represents the noise, accounting for uncertainties in the lever-arm model and in the estimation of the torques from the motion capture data. 

The force exerted by a muscle depends not only on its physiological properties such as mass and length, but also on the state of extension which may depend on time. We refer to the maximum force as the tetanic force. We denote by $f^*_j(t)$ the maximum possible force of the $j$th muscle at time $t$. The muscle force $f_j(t)$ can be expressed relative to the tetanic force, writing
    \begin{align}
     f_j(t) = q_j(t) f^*_j(t), \quad 0\leq q_j(t) \leq 1,
    \end{align}
the function $q_j(t)$ being referred to as the muscle activation.  We express the torques in terms of the activation functions as
    \begin{equation}\label{model}
       \bfb(t) = \bfA(t) \bff(t) + \bfe(t)  = \underbrace{\bfA(t) \bfF^*(t)}_{=\bfA^*(t)} \bfq(t) + \bfe(t) = \bfA^*(t) \bfq(t) + \bfe(t),
     \end{equation}
where $\bfF^*(t) \in\Rbb^{n_m\times n_m}$ is a diagonal matrix with the tetanic forces on its diagonal. The {\em muscle recruitment problem} is to characterize all activation configurations $\bfq(t)$ satisfying  equation \eqref{model}.

\section{Details of the Myobolica Prior}
\label{sec:myobolica_prior}
The Myobolica prior, $\pi_{\rm prior}$, encodes physiologically relevant properties of the muscle activations. Given the time-discretized activation matrix $\bfQ$  (\ref{Q matrix}), we write the prior model as 
	\begin{align}
		\pi_{\rm prior}(\bfQ = \pi_{\rm box}(\bfQ) \pi_{\rm min}(\bfQ) \pi_{\rm yank}(\bfQ),
	\end{align}
where the first factor,  $\pi_{\rm box}$, enforces the torques to obey the physiologically relevant bound constraints, that is, 
	\begin{align}
		\pi_{\rm box}(\bfQ) 
			= \begin{cases} 
				1, & \text{if $ \bfq(t_k) \in [0,1]^\nmuscle$ $1\leq k\leq n_t$,}\\
				0, & \text{otherwise}
			\end{cases}
	\end{align}
where $[0,1]^\nmuscle$ denotes the unit hypercube in $\Rbb^{\nmuscle}$. 
The second factor, $\pi_{\rm min}$, favors muscle activations with lower activity levels, thereby encouraging the more energy efficient walking patterns, 
	\begin{align}
	\pi_{\rm min}(\bfQ) \propto \prod_{k=1}^{n_t}\exp\left( -\tfrac{1}{2\alpha^2} \|\bfq(t_k)\|_2^2\right)
	\end{align}
where $\alpha > 0$ is a tunable parameter controlling the bias toward low activity. Observe that these factors treat the time slices as a priori independent random variances. The temporal correlation is introduced by the third factor, $\pi_{\rm yank}$, encourages smooth transitions of muscle activations over time. Denoting by $\bfF_j$ the $n_t$-vector with components $F_j(t_k) = F_j^* q_j(t_k)$, $1\leq k\leq n_t$,
we write the prior as
	\begin{align}
		\pi_{\rm yank}(\bfQ) \propto \prod_{j=1}^{n_m}\exp\left(- \tfrac{1}{2\beta^2} \|\bfL \bfF_j\|_2^2\right),
	\end{align}
where $\beta>0$ is a tuning parameter,  and $\bfL$ is a first order finite difference operator with appropriate boundary conditions imbedded in it. In \cite{amankwah_exploring_2024}, the effect of the boundary conditions (fixed, random, or periodic) are discussed in more detail.  


\section{CP Factorization Algorithms}
\label{sec:cp_algorithms}
We present vanilla CP and NNCP factorization algorithms. 
For more details, see \cite{ballard_tensor_2025, kolda_tensor_2009, bader_matlab_2017}

\begin{algorithm}
\caption{Nonnegative CP Decomposition via Multiplicative Updating}\label{alg:cp_nmu}
    \begin{algorithmic}[1]
        \State {\bf Inputs:} nonnegative muscle activation tensor $\TQ\in \Rbb_+^{\nmuscle \times \ntime\times \nsample}$, user-prescribed rank $R_{\rm syn}$
        \State {\bf Return:} CP factor matrices $\Umuscle\in \Rbb_+^{\nmuscle\times R_{\rm syn}}$, $\Utime\in \Rbb_+^{\ntime\times R_{\rm syn}}$, $\Usample\in \Rbb_+^{\nsample\times R_{\rm syn}}$ 
        \Statex
        \State Initialize $\Umuscle^{(0)}$, $\Utime^{(0)}$, $\Usample^{(0)}$ with nonnegative entries
        \For{$i=0,1,2, \dots$}  
            \State Construct helper matrix $\bfH_m = \Usample^{(i)} \odot \Utime^{(i)}$ and update 
            	$$\Umuscle^{(i+1)} = \Umuscle^{(i)}* \left(\frac{\bfQ_{(1)} \bfH_m}{ \Umuscle^{(i)}\bfH_m^\top \bfH_m + \varepsilon}\right)$$
            \State Construct helper matrix $\bfH_t = \Usample^{(i)} \odot \Umuscle^{(i + 1)}$ and update 
            	$$\Utime^{(i+1)} = \Utime^{(i)}* \left(\frac{\bfQ_{(2)} \bfH_t}{\Utime^{(i)}\bfH_t^\top \bfH_t + \varepsilon}\right)$$
            \State Construct helper matrix $\bfH_s = \Utime^{(i + 1)} \odot \Umuscle^{(i + 1)}$ and update
            	 $$\Usample^{(i+1)} = \Usample^{(i)}* \left(\frac{\bfQ_{(3)} \bfH_s}{ \Usample^{(i)}\bfH_s^\top \bfH_s + \varepsilon}\right)$$
            \Statex
        \EndFor
    \end{algorithmic}
\end{algorithm}

Note that in practice, we postprocess the CP factor matrices by normalizing and permuting the columns of the NNCP factor matrices based on order of magnitude (i.e., the first column of each NNCP factor matrices is associated with the largest normalization constant and so on).

\begin{algorithm}
\caption{CP Decomposition via Alternating Least Squares}\label{alg:cp_als}
    \begin{algorithmic}[1]
        \State {\bf Inputs:} Bayesian approximation error tensor $\TE\in \Rbb_+^{\nmuscle \times \ntime\times \nsample}$, user-prescribed rank $R_{\rm err}$
        \State {\bf Return:} CP factor matrices $\bfU_{m}\in \Rbb_+^{\nmuscle\times R_{\rm err}}$, $\bfU_{t}\in \Rbb_+^{\ntime\times R_{\rm err}}$, $\bfU_{s}\in \Rbb_+^{\nsample\times R_{\rm err}}$ 
        \Statex
        \State Initialize $\bfU_m^{(0)}$, $\bfU_t^{(0)}$, $\bfU_s^{(0)}$  
        \For{$i=0,1,2, \dots$}  
            \State Solve mode-$1$ unfolded least-squares problem 
            	$$\bfU_m^{(i+1)} = \operatorname*{argmin}_{\bfX_m \in \Rbb^{\nmuscle\times R_{\rm err}}} \|\bfE_{(1)} - \bfX_m (\bfU_s^{(i)} \odot \bfU_t^{(i)})^\top\|_F^2$$
	     \State Solve mode-$2$ unfolded least-squares problem 
            	$$\bfU_t^{(i+1)} = \operatorname*{argmin}_{\bfX_t \in \Rbb^{\nmuscle\times R_{\rm err}}} \|\bfE_{(2)} - \bfX_t (\bfU_s^{(i)} \odot \bfU_m^{(i+1)})^\top\|_F^2$$
	\State Solve mode-$3$ unfolded least-squares problem 
            	$$\bfU_s^{(i+1)} = \operatorname*{argmin}_{\bfX_s \in \Rbb^{\nmuscle\times R_{\rm err}}} \|\bfE_{(3)} - \bfX_s (\bfU_t^{(i+1)} \odot \bfU_m^{(i+1)})^\top\|_F^2$$
	\EndFor
    \end{algorithmic}
\end{algorithm}

\section{Linear Discriminant Analysis}
\label{sec:lda_details}
In this section, be briefly outline the LDA method of finding a projection direction that optimally separates the projections of two samples. Consider two multidimensional random samples
    \begin{align*}
     {\Scal} =\{\bfx^{(1)}, \ldots, \bfx^{(N)}\},\quad {\Scal}' =\{\bfy^{(1)}, \ldots, \bfy^{(M)}\},
    \end{align*}
where $\bfx^{(j)}, \bfy^{(j)}\in{\mathbb R}^n$ are drawn from two possibly different probability distributions, and the goal is to assess, based on the given samples, how close the underlying distributions are.  In the present context, the sample ${\Scal}$ represents the Myobolica sample of activations of a selected muscle, while ${\Scal}'$ comprises the corresponding sample generated by the tensor-based emulator.
Given the sample means,
    \begin{align*}
     \overline \bfx = \frac 1N\sum_{j=1}^N \bfx^{(j)}, \quad \overline \bfy = \frac 1M\sum_{j=1}^M \bfy^{(j)}, \quad
     \overline \bfc = \frac 1{N+M}\left(\sum_{j=1}^N \bfx^{(j)} +
     \sum_{j=1}^M \bfy^{(j)}\right),
    \end{align*}
we define two auxiliary samples of size $N+M$, the first one by replacing each $\bfx^{(j)}$ and $\bfy^{(j)}$ by their respective sample mean,
    \begin{align*}
     \bfc^{(j)} = \left\{\begin{array}{ll}
    \overline \bfx-\overline \bfc, & 1\leq j\leq N,\\
    \overline \bfy-\overline \bfc, & N+1\leq j\leq N+M,
     \end{array}
     \right.
    \end{align*}
and the second one obtained by centering the samples ${\Scal}$ and ${\Scal}'$,
    \begin{align*}
     \bfz^{(j)} = \left\{\begin{array}{ll}
    \bfx^{(j)} - \overline \bfx, & 1\leq j\leq N,\\
    \bfy^{(j-N)} - \overline \bfy, & N+1\leq j\leq N+M,
     \end{array}
     \right.
    \end{align*}

We denote by ${\bfC}\in{\mathbb R}^{n\times(N+M)}$ and
${\bfZ}\in{\mathbb R}^{n\times(N+M)}$ the matrices with columns $\bfc^{(j)}$ and $\bfz^{(j)}$, respectively. It is assumed that the samples are extensive enough to guarantee that the matrix ${\bfZ}$ is of full rank $n$. Given a projection direction ${\bfxi}\in{\mathbb R}^n$, in the LDA analysis the goal is to choose the direction in such a way that the univariate sample $\{{\bfxi}^\top \bfc^{(j)}\}$ has the largest possible variance, while at the same time, the sample $\{{\bfxi}^\top \bfz^{(j)}\}$ has the smallest possible variance. This multi-objective optimization problem is replaced by the problem of maximizing the objective function
\[
 H({\bfxi}) = \frac{{\rm var}\{{\bfxi}^\top \bfc^{(j)}\}} {{\rm var}\{{\bfxi}^\top \bfz^{(j)}\}} 
 = \frac
 {{\bfxi}^\top{\bfC} {\bfxi}}
 {{\bfxi}^\top{\bfZ} {\bfxi}}.
\]
The idea is to find the projection direction in which the two subsamples ${\Scal}$  and ${\Scal}'$ appear as compact as possible (minimum within-class variance of 
$\{{\bfxi}^\top \bfz^{(j)}\}$), while being separated from each other as well as possible (maximum between-class variance of 
$\{{\bfxi}^\top \bfc^{(j)}\}$). It turns out that the optimal direction ${\bfxi}$ can be found as an eigenvector corresponding to the largest generalized eigenvalue $\lambda$ satisfying
\[
 {\bfC}{\bfxi} = \lambda{\bfZ}{\bfxi}.
\]
Observe that the process described here is purely data-driven, and no underlying assumptions about the probability distributions are made.
The generalized eigenvalue problem is equivalent to the standard eigenvalue problem
\[
 {\bfK}^{-\top}{\bfC}{\bfK}^{-1} {\bfw} = \lambda {\bfw}, \quad {\bfw} = {\bfK}{\bfxi},
\]
where ${\bfK}$ is the Cholesky factor of the matrix ${\bfZ}^\top{\bfZ}$. 
For more details, we refer to \cite{calvetti_mathematics_2020}. Once the optimal projection direction is determined, the comparison of the samples ${\Scal}$ and ${\Scal}'$ is performed by considering the univariate samples
    \begin{align*}
     {\Scal}_{\bfxi} =\{{\bfxi}^\top\bfx^{(1)}, \ldots, {\bfxi}^\top\bfx^{(N)}\},\quad {\Scal}'_{\bfxi} =\{{\bfxi}^\top\bfy^{(1)}, \ldots, {\bfxi}^\top\bfy^{(M)}\},
    \end{align*}
e.g., by considering the histograms.

\printbibliography

\end{document}